\documentclass{article}

\usepackage{amsfonts, amssymb, palatino, amscd}

\newtheorem{theorem}{Theorem}
\newtheorem{lemma}[theorem]{Lemma}

\newtheorem{prop}[theorem]{Proposition}
\newtheorem{corr}[theorem]{Corollary}
\newtheorem{conj}[theorem]{Conjecture}

\newcommand{\cyc}{\langle g \rangle}
\begin{document}

\title{Orbifold Cohomology as Periodic Cyclic Homology}
\author{Vladimir Baranovsky}

\date{June 22, 2002}

\maketitle

\begin{abstract}
It known from the work of Feigin-Tsygan, Weibel and Keller that the 
cohomology groups of a smooth complex variety $X$ can be recovered
from (roughly speaking) its derived category of coherent
sheaves. In this paper we show that for a finite group $G$ acting
on $X$ the same procedure applied to $G$-equivariant sheaves
gives the orbifold cohomology of $X/G$.

As an application, in some cases we are able to obtain simple
proofs of an additive isomorphism  between the
orbifold cohomology of $X/G$ and the usual cohomology of its
crepant resolution (the equality of Euler and Hodge numbers
was obtained earlier by various authors). We also state some
conjectures on the product
structures, as well as the singular case; and a connection with 
a recent work by Kawamata.

\end{abstract}

\section{Introduction}

Let $X$ be a smooth variety over a field $k$ (for simplicity we
assume in this introduction that $k = \mathbb{C}$) and $G$ be a
finite group acting on $X$. If the quotient variety $X/G$ is Gorenstein
(i.e. the canonical class is a Cartier divisor) and $\pi: Y \to X/G$
is a crepant resolution of singularities then Ruan's Cohomological
Crepant Resolution Conjecture (which we call Cohomological Conjecture
for short) states, as a particular case, that the cohomology groups
$H^*(Y, \mathbb{C})$ should be isomorphic to the orbifold cohomology
$$
H^*_{orb}(X/G, \mathbb{C}) = \Big( \bigoplus_{g \in G} H^*(X^g, \mathbb{C})
\Big)_G
$$
where $(\ldots)_G$ denotes the coinvariants, $X^g$ is the fixed
point set, and $G$ acts on the above direct sum by conjugating $g$.
This definition of $H^*_{orb}(X, \mathbb{C})$ is slightly different
from the usual one, see Section 3 of \cite{Ru}, but equivalent to it.
Moreover, Ruan has introduced product structures on
$H^*(Y, \mathbb{C})$ (a deformation of the usual product using the
rational
curves contracted by $\pi$) and $H^*_{orb}(X/G, \mathbb{C})$,
\textit{loc. cit.},
and the Cohomological Conjecture states that $H^*(Y, \mathbb{C})
\simeq H^*_{orb}(X/G, \mathbb{C})$ is actually a ring isomorphism.
On the level of Betti (or Hodge) numbers this conjecture was recently
proved by Lupercio-Poddar and Yasuda, see \cite{Y}. However, with the
approach
used in the proof (motivic integration) it is not clear how to
identify the actual cohomology groups with their product structures.
In this paper we will try to outline a different approach to the
Cohomological
Crepant Resolution Conjecture and show that it is in fact a consequence
of a Categorical Resolution Conjecture stated (in a form and under
a name slightly different from ours) by Kawamata in \cite{Ka}. We hope
that the categorical approach will allow to interpret
the product structures. Besides, the author believes that it would
be very important to establish a link between the
categorical and the motivic integration
approach; relating the derived category of sheaves to the space
of arcs (or possibly the space of formal loops of Kapranov-Vasserot).

For simplicity we only work with global quotients $X/G$ but
all statements can be made (and, hopefully, proved)
for general smooth Deligne-Mumford stacks and categories of
sheaves twisted by a gerbe.

In more detail: we would like to show that the above Cohomological Conjecture
follows from an equivalence of two derived categories: the bounded
derived category $D^b(Y)$ of coherent sheaves  on $Y$ and the
bounded derived category $D^b_G(X)$ of $G$-equivariant sheaves on $X$
(i.e. sheaves on the quotient stack $[X/G]$).
Thus, a possible proof of the Cohomological Conjecture could consist
of three steps:

(1) Prove an equivalence of derived categories $D^b(Y) \to D^b_G(X)$
(Categorical Resolution Conjecture - see end of Section 5).

(2) Recover an isomorphism of cohomology groups
$H^*(Y, \mathbb{C}) \to H^*_{orb}(X/G, \mathbb{C})$ from the above
equivalence.

(3) Identify the two product structures.

\medskip

\noindent
In this paper we mostly deal with the second step. As for the
other two, we note that (1) is known in some cases
due to the work of Bridgeland-King-Reid and
Kawamata (see \cite{BKR}, \cite{Ka} and Section 4 of this paper);
while the orbifold product in (3)
seems to arise from the convolution product of sheaves
(see Section 5).

To deal with (2) one needs a construction
which recovers the (orbifold) cohomology ring from the (equivariant) 
derived category. In a sense, we are using some additional structure:
the derived category should "remember" that it was obtained
as a quotient of two DG-categories forming a \textit{localization pair},
see Section 2.4 of \cite{K1}. However, it follows from a result of Orlov
\cite{Or} that fully faithful exact functor between the (non-equivariant)
derived
categories automatically preserves this additional structure; and
the same holds for equivalences of equivariant
derived categories (in the case of a finite group action).

In the non-equivariant case a construction recovering the cohomology groups
follows from the work of Feigin-Tsygan, Weibel and Keller. In fact,
for any exact category $\mathcal{A}$, such as the category $Vect(Y)$ of vector
bundles on $Y$ or the category $Vect_G(X)$ of $G$-equivariant vector bundles on
$X$, Keller constructs in \cite{K1} a \textit{mixed complex} $C(\mathcal{A})$ which
leads to a family of homology theories $HC_{\bullet}(\mathcal{A}, W)$
depending on a graded $k[u]$-module $W$ (usually taken to be of finite
projective dimension). We "recall" the relevant definitions in Section 2.
Two properties make this construction very attractive in our setup:
\begin{itemize}

\item  When $\mathcal{A}$ is the category of vector bundles on $Y$, $k = 
\mathbb{C}$  and  $W=k[u, u^{-1}]$, the above homology group 
$HC_0(\mathcal{A}, W)$ (resp. $HC_1(\mathcal{A}, W)$) can be identified
with $H^{even}(Y, \mathbb{C})$ (resp. $H^{odd}(Y, \mathbb{C})$). Note that
multiplication by $u$ gives an isomorphism $HC_i(\mathcal{A}, W) \simeq
HC_{i+2}(\mathcal{A}, W)$.
 
\item For any $W$ the homology theory $HC_{\bullet}(\mathcal{A}, W)$ is
invariant with respect to equivalences of derived categories
coming from functors between localization pairs.
\end{itemize}

\bigskip
\noindent
To formulate our results, recall that $G$ acts on $\coprod_{g \in G} X^g$:
$h \in G$ sends $x \in X^g$ to $hx \in X^{hgh^{-1}}$. 
This action is inherited by  $\bigoplus_{g \in G} HC_{\bullet} 
(Vect(X^g), W)$.

\begin{theorem} \label{main} Let $G$ be a finite group acting on a smooth 
quasiprojective variety $X$ over a field $k$ of characteristic not 
dividing 
$|G|$. For any graded $k[u]$-module $W$ of finite projective dimension 
there exists an isomorphism functorial with respect to pullbacks
under $G$-equivariant maps:
$$
\psi_X: HC_{\bullet} (Vect_G(X), W) \simeq
\Big(\bigoplus_{g \in G} HC_{\bullet}(Vect(X^g), W)\Big)_G
$$
where $(\ldots)_G$ denotes the coinvariants.
\end{theorem}

In the $C^{\infty}$-manifold or $C^{\infty}$-etale groupoid setting
this result (formulated in terms of modules over smooth functions rather
than categories) has a long history. First Feigin-Tsygan, see \cite{FT},
constructed a spectral sequence computing cyclic homology of a general crossed
product algebra, which was later reformulated by Getzler-Jones, see
\cite{GJ2}; and also \cite{AK} for more general crossed products by Hopf
algebras.  When the crossed product algebra comes from functions
on a smooth manifold the $E_2$ term of this spectral sequence can be
interpreted in terms of fixed point submanifolds: this result was
announced in \cite{Br} and the first published proof appears in \cite{BN}.
Later it was generalized in \cite{Cr}.
The case when $G$ is a Lie group was studied by Nistor in \cite{N}; later
Block and Getzler have related the corresponding crossed product cyclic
homology groups to equivariant differential forms and fixed points, see
\cite{BG}. We also mention a
closely related computation of $G$-equivariant topological $K$-theory
by G. Segal, see \cite{HH}; and its algebraic $K$-theory counterparts
\cite{V}, \cite{To}.

In Section 3 we adapt the proofs in \cite{BG} and \cite{GJ2}
to fit our case of categories and rings of regular functions. Note that
the proof of \cite{V} cannot be applied in our case due to the failure of
devissage, see Example 1.11 in \cite{K1}.

In those cases when the derived equivalence $D^b(Y) \to D^b_G(X)$  is known,
we get an isomorphism between $HC_\bullet(Vect(Y), W)$ and the right
hand side expression in Theorem \ref{main}, obtaining a
slightly generalized version of the Cohomological Conjecture (in general
the cyclic homology groups \textit{do not} satisfy the long exact
sequence hence even
on the level of dimensions the equality cannot be derived using
motivic measures and motivic integration). This is the second main
result of this paper (see Corollary \ref{last}).

We expect that, in order to identify the product structures
in step (3) above, one should
modify $\psi_X$ to make it compatible with pushforwards under
$G$-equivariant closed embeddings, rather than pullbacks
(compare with Lemmas 4.2 and 4.3 in \cite{V}).

We note here that for a general algebraic group $G$ the equivariant derived
category $D^b_G(X)$ is defined \textit{not} by taking complexes
of $G$-equivariant sheaves but by a more delicate localization procedure,
see \cite{BL}. One can expect that the corresponding
cyclic homology groups satisfy nice properties, for example
similar to those proved in \cite{BG}.

This paper is organized as follows. Section 2 gives  some basic
information of cyclic homology groups of exact categories. In Section 3
we prove Theorem \ref{main}. In Section 4 we show how an equivalence
of derived categories implies equality of (orbifold) cohomology groups
and also give some examples in which this equivalence is known.
Finally, in Section 5 we give a conjecture about the singular case
and a conjecture  on how the orbifold cohomology product can be recovered
from the convolution product in the derived category.

\medskip
\noindent
\textbf{Acknowledgements.} The present work was motivated by a lecture
of Y. Ruan on orbifold cohomology given by him at Caltech;
and the two beautiful papers by B. Keller
\cite{K1}, \cite{K2} on  cyclic homology of exact categories. The author is
grateful to both of them for providing this inspiration.

\section{Generalities on Mixed Complexes}

Recall, cf \cite{W1} that a \emph{mixed complex} over a
commutative ring $k$ is a sequence of $k$-modules 
$\{C_m: m \geq 0\}$ with two families of morphisms $b: C_m \to C_{m-1}$
and $B: C_m \to C_{m+1}$ satisfying
$b^2 = B^2 = Bb + bB = 0$. To any such  mixed complex
one can apply the following formalism (see \cite{GJ1}): let $W$ be
a graded module over the polynomial ring $k[u]$, where $\deg(u) = -2$
(in practice it is always assumed that $W$ has finite homological 
dimension). Then one can form a complex $C[[u]] 
\otimes_{k[u]} W$ with a differential $b + u B$ and compute its 
cohomology groups, to be denoted by $HC_{\bullet}(C, W)$. The 
following are important examples:

\begin{itemize}
\item $W = k[u]/uk[u]$ gives the Hochschild homology $HH_{\bullet} 
(C)$

\item $W = k [u, u^{-1}]/u k[u]$ gives cyclic homology $HC_{\bullet} 
(C)$

\item $W = k[u, u^{-1}]$ gives periodic cyclic homology $HP_{\bullet} 
(C)$

\item $W = k[u]$ gives negative cyclic homology $HN_{\bullet} 
(C)$ (sometimes also denoted by $HC^-_{\bullet} (C)$).
\end{itemize}

\noindent
The following lemma shows that for some purposes it suffices
to consider only the first case.

\begin{lemma} \label{mix-der}
Let $f: (C, b, B) \to (C', b', B')$ be a map of mixed
complexes such that $f$ induces an isomorphism $H(C, b) \to H(C', b')$.
Then for any coefficients $W$ of finite projective dimension 
over $k[u]$,
$$
f: \qquad H_{\bullet} (C[[u]] \otimes_{k[u]} W, b + u B) \to 
H_{\bullet}(C'[[u]] \otimes_{k[u]} W, b' + u B')
$$
is an isomorphism.
\end{lemma}
\noindent
\textit{Proof.} See Proposition 2.4 in \cite{GJ1}.

\medskip
\noindent
This lemma justifies the following point of view on mixed complexes, 
see \cite{K1}, Section 1.2. Let $\Lambda$ be the DG-algebra
generated by an
indeterminate $\varepsilon$ of chain degree 1 with $\varepsilon^2 = 0$
and $d \varepsilon = 0$. Then a mixed complex may be identified 
with a left $\Lambda$-module whose underlying DG $k$-module is $(C, b)$
and where $\varepsilon$ acts by $B$. Moreover, if we are interested
only in the resulting homology groups (as is the case in this paper),
we can view a mixed complex as an object in the derived category
of the DG algebra $\Lambda$. 

In what follows we will need a definition of the \textit{mapping cone}
over a  map $f: C \to C'$ of mixed complexes. It is given by the 
mixed complex 
$$
\bigg(C' \oplus C[1], \bigg[ 
\begin{array}{cc} b_{C'} & f \\ 0 & - b_C \end{array}
                 \bigg], \bigg[
\begin{array}{cc} B_{C'} & 0 \\ 0 & - B_{C} \end{array}
                 \bigg] \bigg)
$$

Now we briefly recall Keller's construction of the mixed complex 
$C(\mathcal{A})$ of an exact category $\mathcal{A}$ over a field $k$ 
(actually in \cite{K1} the complex is defined for any commutative
ring $k$ but the general definition is somewhat more involved).
Starting from $\mathcal{A}$ one can construct its category
$\mathcal{C}^b \mathcal{A}$ of all bounded complexes over $\mathcal{A}$
and the category $\mathcal{A}c^b \mathcal{A}$ of bounded acyclic
complexes over $\mathcal{A}$. Both $\mathcal{C}^b\mathcal{A}$ and 
$\mathcal{A}c^b \mathcal{A}$ are  \textit{DG categories}, i.e. for any 
pair of objects $X$, $Y$ the group $Hom(X, Y)$ is $\mathbb{Z}$-graded 
(by degree of a map) with a differential, which satisfies some natural 
axioms (see \cite{K3})). 

For any small DG category $\mathcal{B}$ over a field $k$ Keller
constructs a mixed complex as follows. Denote for
notational convenience $Hom_{\mathcal{B}}(X, Y)$ by 
$(X \to Y)$ and consider for each $n \in \mathbb{N}$ a vector space
$$
C_n (\mathcal{B})= \bigoplus (B_0 \to B_1) \otimes (B_1 \to B_2)
\otimes \ldots \otimes (B_{n-1} \to B_{n}) \otimes (B_n \to B_0)
$$
where the sum runs over all sequences $B_0, \ldots, B_n$ of objects
of $\mathcal{B}$. The face maps 
$$
d_i (f_0, \ldots, f_i, f_{i+1}, \ldots, f_n) = \bigg\{
\begin{array}{lll}
(f_0, \ldots, f_{i} f_{i+1}, \ldots, f_n) & \textrm{if } 0 \leq i \leq n-1 \\
(-1)^{n+ \sigma} (f_n f_0, \ldots, f_{n-1}) & \textrm{if } i = n
\end{array}
$$
(where $\sigma = \deg f_n \cdot (\deg f_0 + \ldots + \deg f_{n-1})$);
together with the degeneracy maps
$$
s_i (f_0, \ldots, f_i, f_{i+1}, \ldots, f_n) =
(f_0, \ldots, f_{i}, id_{B_i},  f_{i+1}, \ldots, f_n) \qquad i = 0, \ldots, n
$$
and the cyclic operator
$$
t_n (f_0, \ldots, f_n) = (-1)^{n + \sigma} (f_n, f_0, f_1,
\ldots, f_{n-1}),
$$
define a mixed complex $(C(\mathcal{B}), b, (1 - t) s N)$ as in \cite{GJ2},
Section 2. Note that, unlike in \cite{K1}, we write $fg$ for a
composition of $f:A \to B$ and $g:B \to C$ instead of the usual 
$gf$. This non-traditional notation will allow us later to match
our formulas with those of Getzler-Jones in \cite{GJ2}. 

Next, for any DG subcategory $\mathcal{C} \subset \mathcal{B}$ one has
a mixed complex
$$
C(\mathcal{C}, \mathcal{B}) = Cone (C(\mathcal{C}) \to C(\mathcal{B})).
$$
We will always use this definition when $(\mathcal{C}, \mathcal{B})$
is a \emph{localization pair} in the sense of \cite{K1}, Section 2.4.
Every such localization pair leads to a triangulated
category $\mathcal{T}= \mathcal{B}/\mathcal{C}$ associated to it.
For example, the pair $(\mathcal{A}c^b \mathcal{A}, \mathcal{C}^b \mathcal{A})$
gives rise to the derived category $\mathcal{T}$ of $\mathcal{A}$. One of
the main properties of the mixed  complex $C(\mathcal{C}, \mathcal{B})$
is expressed in the following proposition (see Theorem 2.4(a) in \cite{K1}):
\begin{prop}
\label{invariance}
Let $(\mathcal{C}, \mathcal{B})$ and $(\mathcal{C}', \mathcal{B}')$ be two
localization pairs and $\mathcal{T}$, $\mathcal{T}'$ their derived
categories.
If $F: \mathcal{B} \to \mathcal{B}'$ is an exact
functor which takes $\mathcal{C}$ to $\mathcal{C}'$ and induces and
equivalence up to factors $\mathcal{T} \to \mathcal{T}'$ (cf. 1.5 in
\cite{K1}), then $F$ induces an isomorphism $C(\mathcal{C}, \mathcal{B}) \to
C(\mathcal{C}', \mathcal{B}')$ in the derived category of $\Lambda$.
\end{prop}

If $\mathcal{A}$ is an exact category then its mixed
complex is defined as
$C(\mathcal{A}c^b \mathcal{A}, \mathcal{C}^b \mathcal{A})$.
We will write $\mathcal{A}c^b(X)$ (resp. $\mathcal{C}^b(X)$) instead
of $\mathcal{A}c^b Vect(X)$ (resp. $\mathcal{C}^b Vect(X))$ and
similarly for $\mathcal{A}c^b_G(X)$ and $\mathcal{C}^b_G(X)$.
Define $C(X) = C(\mathcal{A}c^b(X), \mathcal{C}^b(X))$; $C_G(X) =
C(\mathcal{A}c^b_G(X), \mathcal{C}^b_G(X))$.

\section{Mixed Complex for a Finite Group Action}

In this section we will prove Theorem \ref{main}. In view of Lemma
\ref{mix-der} and the subsequent remarks, it suffices to prove
the following proposition:
\begin{prop}
\label{actual}
Let $G$ be a finite group acting on a smooth 
quasiprojective variety $X$ over a field $k$ of characteristic not dividing 
$|G|$. There exists a quasiisomorphism in the derived category of $\Lambda$:
$$
\psi_X: C_G(X) \to
\Big(\bigoplus_{g \in G} C(X^g)\Big)_G
$$
which is functorial with respect to pullbacks under $G$-equivariant
maps of smooth varieties $f: Y \to X$.
\end{prop}

\noindent
The proof of the above proposition will be carried out in several 
steps: first we replace $Vect_G(X)$ by a categorical analogue of
a crossed product ring $A \rtimes G$ and construct a functorial map
of objects in the derived category (Steps 1 and 2). Once this map is
constructed, we use Mayer-Vietoris, Luna's Fundamental Lemma and etale
descent of Weibel-Geller to reduce to the case when $X$ is a vector
space with a linear action of $G$; and conclude the argument using
exactness of a Koszul complex (Steps 3-6). Our proof is an adaptation of
the argument in \cite{BG} from the differentiable to the algebraic situation.

\vspace{1cm}
\noindent
\label{step1}
\textit{Step 1: Change of categories.}

\vspace{1cm}
\noindent
Note that the action of $G$ on $X$ is inherited by the categories
$Vect_G(X)$, $\mathcal{C}^b(X)$ and $\mathcal{A}c^b(X)$:
for any $g \in G$ and any bundle (or a complex of bundles) $\mathcal{F}$
we can consider $g\mathcal{F} := (g^{-1})^* \mathcal{F}$ and
for any $\psi: \mathcal{F} \to \mathcal{G}$ we have $g(\psi): g\mathcal{F} \to g\mathcal{G}$.

For any $\mathcal{F}$  consider
the object $\widetilde{\mathcal{F}} = \bigoplus_{g \in G} g\mathcal{F}$
with its natural $G$-equivariant structure. Then  each object
$\mathcal{H}$ in $\mathcal{C}^b \;Vect_G(X)$
is isomorphic to a direct factor of some $\widetilde{\mathcal{F}}$. 
In fact, take $\mathcal{F} = \mathcal{H}$ viewed as an object in 
$\mathcal{C}^b\;Vect(X)$. Then the $G$-equivariant structure on $\mathcal{H}$
defines isomorphisms $i_g: \mathcal{H} \to g\mathcal{H}$ for
all $g \in G$. Now consider the $G$-equivariant maps
$$
\mathcal{H} \stackrel{a}\longrightarrow \widetilde{\mathcal{H}} = 
\bigoplus_{g \in G} g\mathcal{H} \stackrel{b}\longrightarrow \mathcal{H}
$$
where $a$ is given by a direct sum of $i_g$ and $b = \frac{1}{|G|} \sum i_g$.
Since $b \circ a = id_{\mathcal{H}}$, $\mathcal{H}$ splits off as
a direct factor of $\widetilde{\mathcal{H}}$.

Denote by $\mathcal{C}^b(X) \rtimes G$, resp. 
$\mathcal{A}c^b(X)\rtimes G$ the full subcategory of 
$\mathcal{C}^b_G(X)$, resp. $\mathcal{A}c^b_G(X)$,
formed by all  $\widetilde{\mathcal{F}}$. Since  
$\mathcal{C}^b(X) \rtimes G$ and 
$\mathcal{A}c^b(X)\rtimes G$ are closed under degree-wise
split extensions and shifts, both are exact DG categories
(by Example 2.2 (a) in \cite{K1}). The above argument shows that
the natural embedding $\mathcal{C}^b(X)\rtimes G \to 
\mathcal{C}^b_G(X)$ induces an equivalence \textit{up to
factors} of  associated derived categories. Now Theorem 2.4 (a)
in \cite{K1} gives
\begin{prop}   \label{crossed}
There exists an isomorphism in the mixed derived category 
$$
C(\mathcal{A}c^b(X)\rtimes G, \mathcal{C}^b(X)\rtimes G) \to C_G(X)
$$
which is functorial with respect to $G$-equivariant pullbacks
of vector bundles.
\end{prop}

\vspace{1cm}
\noindent
\label{step2}
\textit{Step 2:  Construction of the Quasiisomorphism.}

\vspace{1cm}
\noindent
To proceed further, we take a closer look at $\mathcal{C}^b(X) 
\rtimes G$. The objects in this category can be identified with 
objects in 
$\mathcal{C}^b(Vect(X))$ while the morphisms are given by 
$$
(\mathcal{F} \to \mathcal{G})_{\mathcal{C}^b(X)\rtimes G}
= \bigoplus_{g \in G} (\mathcal{F} \to 
g\mathcal{G})_{\mathcal{C}^b(X)}
$$
(the expression on the right obviously coincides with all
$G$-equivariant morphisms from $\widetilde{\mathcal{F}}$ to 
$\widetilde{\mathcal{G}}$). If we denote by $\psi \cdot g$ and
$ \varphi \cdot h $
the components of 
$(\mathcal{F} \to \mathcal{G})_{\mathcal{C}^b(X)\rtimes G}$ and
$(\mathcal{G} \to \mathcal{H})_{\mathcal{C}^b(X)\rtimes G}$
living in 
$(\mathcal{F} \to g\mathcal{G})_{\mathcal{C}^b(X)}$ and
$(\mathcal{G} \to h\mathcal{H})_{\mathcal{C}^b(X)}$,
respectively, then the composition of $\psi \cdot g$ and
$\varphi \cdot h $ is given by
$$
\psi g(\varphi) \cdot g h:  \mathcal{F}  \to 
 g\mathcal{G} \to gh\mathcal{H}
$$
Now define the map
$$
C(\mathcal{C}^b(X) \rtimes G) \to \Big(\bigoplus_{g \in G}
C(\mathcal{C}^b(X^g)) \Big)_G
$$
by sending $(\varphi_0 \cdot g_0, \varphi_1 \cdot g_1,
\ldots, \varphi_n \cdot g_n)$ to
$({\varphi_0}_{|_{X^g}},  {g_0(\varphi_1)}_{|_{X^g}},
g_0 g_1(\varphi_2)_{|_{X^g}},  \ldots,  (g_0 \ldots g_{n-1})
(\varphi_n)_{|_{X^g}})$ where $g = g_0 \ldots g_n$. Note that only
the map to the coinvariants is indeed a map of mixed complexes
since the individual maps $C(\mathcal{C}^b(X) \rtimes G) \to
C(\mathcal{C}^b(X^g))$ do not respect the last face map $d_n$ and
the cyclic  operator $t_n$. Since the above morphism obviously sends
$
C(\mathcal{A}c^b(X) \rtimes G)$ to $\Big(\bigoplus_{g \in G}
C(\mathcal{A}c^b(X^g)) \Big)_G
$
we obtain a morphism of mixed complexes
$$
C(\mathcal{A}c^b(X) \rtimes G, \mathcal{C}^b(X) \rtimes G) \to
\Big(\bigoplus_{g \in G} C(\mathcal{A}c^b(X^g),
\mathcal{C}^b(X^g)) \Big)_G.
$$
Composing this map with the inverse of the quasiisomorphism
in Proposition \ref{crossed} we obtain a morphism of objects
in the derived category of $\Lambda$
$$
\psi_X: C_G(X) \to \Big( \bigoplus_{g \in G} C(X^g) \Big)_G
$$
which is functorial with respect to pullbacks under $G$-equivariant
maps $f: Y \to X$.

\vspace{1cm}
\noindent
\label{step3}
\textit{Step 3: Mayer-Vietoris and Luna's Fundamental Lemma}

\vspace{1cm}
\noindent
Now we use a Mayer-Vietoris sequence to reduce to the case when
$X$ is affine.

\begin{prop}   Let $X$ be a quasiprojective scheme,
$V, W \subset X$ two $G$-invariant open subschemes and
 $U = V \cap W$. There is a distinguished triangle in the
mixed derived category category:
$$
C_G(X) \to C_G(V) \oplus C_G(W) \to C_G(U) \to C_G(X)[1]
$$
\end{prop}
\noindent
\textit{Outline of Proof.} Most of the argument is identical to
the non-equivariant case proved in  Proposition 5.8 of \cite{K2}.
First one shows that for a quasiprojective scheme $X$,
$C_G(X)$ is quasiisomorphic to the complex
obtained from the category of $G$-equivariant  perfect complexes
(see Section 5.1 of \cite{K2}). Moreover, let $\mathcal{T}_G(X)$
be the derived category of $G$-equivariant perfect complexes and for
any closed $Z \subset X$ by $\mathcal{T}_G(X \textrm{on } Z)$ the
subcategory of complexes which are exact on the complement to $Z$.
If $Z = X \setminus W$
and $j: V \to X$ is the open embedding, one shows that the lines
of the diagram
$$
\begin{CD}
0 @>>> \mathcal{T}_G (X \textrm{ on } Z) 
  @>>> \mathcal{T}_G  X @>>> \mathcal{T}_G W @>>> 0\\
& & @Vj^*VV @VVV @VVV  \\
0 @>>> \mathcal{T}_G (V \textrm{ on } Z) 
  @>>> \mathcal{T}_G V @>>> \mathcal{T}_G (V \cap W) @>>> 0\\
\end{CD}
$$
are exact up to factors and the functor $j^*$ is an 
equivalence up to factors (the proof in Sections 5.4 and 5.5 of
\cite{TT} may be repeated almost word-by-word; note also
that  working ``up to factors" allows one to ignore all problems
with non-surjectivity of $K^0_G(X) \to  K^0_G(W)$  and so on,
since for any complex $E$, the sheaf $E \oplus E[1]$ has
zero class in $K$-theory). Finally, one applies Theorem 2.7
of \cite{K1}.

Alternatively, one can work with $\mathcal{C}^b(X) \rtimes G$ and
deduce the above fact from the non-equivariant version (\cite{K1},
Proposition 5.8)
$\square$

The above Mayer-Vietoris sequence shows that if $\psi_U$, $\psi_V$
and $\psi_W$ are quasiisomorphisms then the same holds for $\psi_X$
(of course, this is only possible by functoriality of $\psi_X$ applied
to open embeddings). By induction on the number of elements in a
$G$-invariant affine covering of $X$ we can assume that $X = Spec\; A$
is affine. In this case we can reinterpret the map $\psi_X$ as
follows. The group $G$ acts on the $k$-algebra $A$ and we can define
$A \rtimes G$ to be the crossed product algebra $A \otimes k[G]$ with
the product defined by
$$
(a_1 \cdot g_1) (a_2 \cdot g_2) = (a_1 g_1(a_2)) \cdot g_1 g_2
$$
Then the category $Vect_G(X)$ is equivalent to the category of
$A \rtimes G$-modules which are projective as $A$-modules. Since $G$
is finite, this is equivalent to
being projective as  $A \rtimes G$-modules.

Now consider the category $dgfree\; A\rtimes G$ of complexes
of free $A\rtimes G$-modules. By Section 2.4 of \cite{K1} the
natural functor $dgfree\; A \rtimes G \to \mathcal{C}^b \; Vect_G(X)$
induces a quasiisomorphism of mixed complexes
$$
C(0, dgfree\; A\rtimes G) \to C_G(X).
$$
Moreover, if one considers $A \rtimes G$ as a subcategory of $dgfree \; A
\rtimes G$ with one object (a free rank one module viewed as
a complex in degree 0) then by Theorem 2.4 (a) of \cite{K1} one gets
a quasiisomorphism
$$
C(A \rtimes G) \to C(0, dgfree \; A \rtimes G)
$$
where the mixed complex $C(B)$ of any $k$-algebra $B$ is defined
as in Section 2 (if one views $B$ as a category with one object).
Thus, we obtain a chain of quasiisomorphisms
$$
C(A \rtimes G) \to C(0, dgfree\; A \rtimes G)
\to C (\mathcal{A}c^b(X) \rtimes G, \mathcal{C}^b (X) \rtimes G)
\to C_G(X)
$$
induced by embeddings of subcategories.

Similarly, if $g \in G$ and $J_g \subset A$ denotes the ideal
of the fixed point set $X^g \subset X$ then we have a chain of
quasiisomorphisms
$$
C(A/J_g) \to C(0, dgfree\; A/J_g) \to C(X^g)
$$
Restricting the map $\psi_X$ constructed in the previous step,
it suffices to prove that the following map is a quasiisomorphism
in the derived category of $\Lambda$:
$$
\psi_A: C(A \rtimes G) \to \Big(\bigoplus_{g \in G} C(A/J_g) \Big)_G;
$$
where $\psi_A(a_0 \cdot g_0, \ldots, a_n \cdot g_n)$ is
the image of $(a_0,  g_0(a_1),  \ldots,
(g_0\ldots g_{n-1})(a_n)) \in A^{\otimes n + 1}$ in
$(A/J_g)^{\otimes n+1} = C_n(A/J_g)$ and $g = g_0 \ldots g_n$.

\bigskip
\noindent
Eventually we will assume that $X$ is not only affine but has some
additional properties. The following proposition is
probably well-known but the author was unable to
find a convenient reference.

\begin{prop}
\label{mayer}
 Let $G$ be a finite group acting on a
smooth quasiprojective variety $X$ over a field $k$. Assume that
$|G|$ is invertible in $k$. There exists a covering of $X$ by
$G$-invariant affine open subsets $U_1, \ldots, U_n$ such that
for any $i = 1, \ldots, n$ there is a point $x_i \in U_i$
satisfying the following properties:

(a)  the fixed point scheme $U_i^g$ is empty unless $g \in G_{x_i}$
(the stabilizer of $x_i$);

(b) if $T_{x_i}$ is the tangent space to $U_i$ at $x_i$ with
its natural $G_{x_i}$ action, then there exists a $G_{x_i}$-equivariant
etale morphism $\varphi_i: U_i \to T_{x_i}$  such that
for any subgroup $H \subset G_{x_i}$ the induced morphism
$U_i/H \to T_{x_i}/H$ is etale and the
diagram
$$
\begin{array}{ccc}
U_i & \longrightarrow & T_{x_i} \\
\downarrow & & \downarrow \\
U_i/H & \longrightarrow & T_{x_i}/H
\end{array}
$$
is Cartesian;

(c) for any $H \subset G_{x_i}$ the fixed point scheme $U_i^H$ is
a scheme-theoretic preimage of $T_{x_i}^H$.
\end{prop}
\noindent
\textit{Proof.} For any point $x \in X$ we can choose and open affine
$G$-invariant neighborhood $U$ such that (a) holds.
By Lemma 8.3 of \cite{BR} there exists a $G_x$-equivariant
map $\varphi: U \to T_x$ such that $\varphi(x) = 0$ and
$d\varphi (x)$ is equal to identity. Now by Theorem 6.2 of \cite{BR}
(Luna's Fundamental Lemma in finite characteristic) applied to
the finite set of subgroups $H \in G_{x_i}$ we can shrink $U$ so
that (b) holds as well.

For (c) let $U_i = Spec\; B$ and $T_{x_i} = Spec\; A$. Denote by $I_A$
the $A^H$-submodule of $A$ generated by elements $a - h(a)$, $a \in A, h\in H$.
Let also $J_A$ be the ideal in $A$ generated by $I_A$; then $J_A$ is the
ideal of the fixed point scheme $T_{x_i}^H$. We use the similar
notation $I_B, J_B$ for the objects corresponding to $B$. Since $B$ is
flat over $A$, $J_A \otimes_A B$ is naturally an ideal in $B = A \otimes_A B$
and it suffices to prove that it coincides with $J_B$ (apriori we
just have an inclusion $J_A \otimes_A B \subset J_B$). Note that since
$|H|$ is invertible $A = I_A \oplus A^H$ as  $A^H$-modules. Applying
$\otimes_{A^H} B^H$ and using the fact that $B = A \otimes_{A^H} B^H$ by
part (b) we conclude that $I_B = I_A \otimes_{A^H} B^H $ which implies
$J_B = J_A \otimes_A B$ by definition of $J_A$, $J_B$. $\square$

\medskip
\noindent
Thus, later we may assume that $X$ is affine and there exists $x \in X$
such that (a), (b) and (c) above are satisfied.

\vspace{1cm}

\noindent
\label{step4}
\textit{Step 4: Eilenberg-Zilber Theorem}

\vspace{1cm}
\noindent
The mixed complex $C(A \rtimes G)$ was studied in a more general
situation by Jones and Getzler \cite{GJ2}. In the next proposition
we present their results (with some simplifications possible due to
the fact that $|G|$ is invertible in $k$; also our isomorphism
is more explicit than in \cite{GJ2}).

To state the lemma we need to fix some notation. For any $g \in G$
consider the sequence of vector spaces
$
\big(A^\natural_g)_n = A^{\otimes n+1}
$
together with the face, degeneracy and cyclic operators defined
by the formulas similar to those in Section 2:
$$
d_i (a_0, \ldots, a_i, a_{i+1}, \ldots, a_n) = \bigg\{
\begin{array}{lll}
(a_0, \ldots, a_{i} a_{i+1}, \ldots, a_n) & \textrm{if } 0 \leq i \leq n-1 \\
(g^{-1}(a_n) a_0, \ldots, a_{n-1}) & \textrm{if } i = n
\end{array}
$$
$$
s_i (a_0, \ldots, a_n) = (a_0, \ldots, a_i, 1, a_{i+1}, \ldots, a_n) \qquad
i = 0, \ldots, n;
$$
$$
t_n(a_0, \ldots, a_n) = (g^{-1}(a_n), a_0, \ldots, a_{n-1})
$$
Note that the operator $t_n$ \emph{does not} satisfy the cyclic
identity $t^{n+1}=1$. However, we can still construct a mixed complex
from the spaces $A^{\natural}_g$ by considering a $G$-action on the
direct sum $\bigoplus_{g \in G} A^\natural_g$ such that $h \in G$
sends $(a_0, \ldots, a_n) \in A^\natural_{g}$ to
$(h(a_0),  \ldots,  h(a_n)) \in A^\natural_{hg h^{-1}}$.
Then this $G$-action commutes with $d_i, s_i, t$ hence we obtain
the face, degeneracy and cyclic operators on the quotient
space $\Big(\bigoplus_g A^\natural_g \Big)_G$ of $G$-coinvariants.
In this quotient space the operator $t_n$ does satisfy $t_n^{n+1} = 1$
and we denote the resulting  mixed complex  by
$C \Big( \bigoplus_{g \in G} A^\natural_g \Big)_G$.

\begin{prop} The map of mixed complexes
$$
 \varphi_A: C(A \rtimes G) \to C\Big( \bigoplus_{g \in G}
 A^\natural_g \Big)_G
$$
defined by $\varphi(a_0 \cdot g_0, \ldots a_n \cdot g_n) =
(a_0, g_0(a_1), g_0 g_1 (a_2), \ldots, g_0 \ldots g_{n-1} (a_n))
\in \big(A^\natural_g\big)_n$ with $g = g_0 \ldots g_n$; is a
quasiisomorphism in the derived category of $\Lambda$.
\end{prop}
\noindent
\textit{Proof.} In \cite{GJ2} Getzler and Jones define
a bi-graded object
$A \natural G (p, q) = k[G^{p+1}] \otimes A^{\otimes q+1}$
with two families of face maps $d^{h}: A \natural G (p, q)
\to A \natural G(p-1, q)$ (horizontal maps) and
$ d^{v}: A \natural G (p, q) \to A \natural G(p, q-1)$ (vertical maps),
and similarly
for degeneracies and cyclic operators (we will not use
the precise definitions of these operators). These two families
of operators give $A\natural G(p, q)$ the structure of a
\emph{cylindrical module}, see \cite{GJ2} before Proposition 1.1.

This cylindrical module has a total complex $Tot_n (A, G) =
\bigoplus_{p + q = n} A \natural G (p, q)$ and a diagonal complex
$\Delta_n(A, G) = A \natural G (n, n)$, both being mixed complexes (see
\cite{GJ2} for more details).
The cylindrical module structure  on $A \natural G(p, q)$ is
defined in such a way that $C(A \rtimes G)$ is isomorphic to
the diagonal complex $\Delta (A, G)$ via the map
$$
(a_0 \cdot g_0, \ldots, a_n \cdot g_n) \mapsto
(g_0, \ldots, g_n | h_0^{-1} a_0, \ldots, h_n^{-1} a_n) \in
A \natural G (n, n) = \Delta_n (A, G)
$$
where $h_i = g_i \ldots g_n$. Applying the Eilenberg-Zilber
Theorem for paracyclic modules, see Theorem 3.1 in \cite{GJ2},
one gets an explicit quasiisomorphism
$$
AW: \Delta (A, G) \to Tot (A, G)
$$
given by the Alexander-Whitney map defined in Section 8.5.2 of
\cite{W1}. We will only need one component of this map
(with values in $A \natural G (0, n)$) which is given simply by
$d_1^h \ldots  d_n^h$, where $d^h_i: A \natural G (i, q)
\to A \natural G (i - 1, q)$ sends $(g_0, \ldots, g_i | a_0, \ldots, a_q)$
to $(g_i g_0, \ldots, g_{i-1} | g_i (a_0), \ldots, g_i (a_0))$.

Now the transformation
$$
(g_0, \ldots, g_p | a_0, \ldots, a_q) \mapsto
(g_1, \ldots, g_p | g_0 g_1 \ldots g_p | a_0, \ldots, a_q)
$$
identifies $A \natural G (p, q)$ with a cylindrical module
the vertical maps of which are given by the above operators on
$\bigoplus_{g \in G} A^\natural_g$
(and the index $g$ corresponds to $(\ldots | g | \ldots )$
in the above notation), while the rows can be identified with the
standard homological complex of  $G$ acting
on $\bigoplus_{g \in G} A^\natural_g$. See Section 4 of \cite{GJ2}
for details.  Since $|G|$ is invertible in $k$, the projection of
$Tot (A, G)$ onto its first column $A \natural (0, \bullet)$
together with the projection to the coinvariants, gives a
quasiisomorphism
$$
Tot(A, G) \to C \Big(\bigoplus_{g \in G} A^\natural_g \Big)_G.
$$
The fact that it is indeed a quasiisomorphism can be proved using
the homotopy
$$
h: A \natural G (p, q) \to A \natural G (p+1, q);
\qquad (g_1, \ldots, g_p | g | a_0 \ldots a_q) \mapsto
\frac{1}{|G|}\sum_{g' \in G} (g', g_1, \ldots, g_p |g
| a_0, \ldots, a_q )
$$
The composition of (quasi)isomorphisms:
$$
C(A \rtimes G) \to \Delta (A, G) \to Tot (A, G) \to
C\Big(\bigoplus_g A^\natural_g \Big)_G
$$
is given by
$$
(a_0 \cdot g_0, \ldots, a_n  \cdot g_n) \to
(g^{-1}_0(a_0), a_1, g_1(a_2), \ldots, (g_1 \ldots g_{n-1}) (a_n))
\in (A^\natural_g)_n; \qquad g = g_1 g_2 \ldots g_n g_0.
$$
Note that this differs from the map in the statement of our
proposition exactly be the action of $g_0 \in G$. Since on the
space of coinvariants $g_0$ acts by identity, this finishes the proof.
$\square$

\bigskip
\noindent
Thus, we have reduced Proposition \ref{actual} to the claim that the
map
$$
C\Big(\bigoplus_{g \in G} A^\natural_g \Big)_G \to
\Big(\bigoplus_{g \in G} C(A/J_g)\Big)_G
$$
defined by the natural surjections $(A^\natural_g)_n  \to (A/J_g)^{\otimes
n+1}$, is a quasiisomorphism. Note also that we have not used the
smoothness assumption yet.

\vspace{1cm}
\noindent
\label{step5}
\textit{Step 5: Shapiro's Lemma and Etale Descent}

\vspace{1cm}
\noindent
Let $\mathcal{O} \subset G$ be a conjugacy class. It is
easy to see that
$
\bigoplus_{g \in \mathcal{O}} A^\natural_g
$
and
$
\bigoplus_{g \in \mathcal{O}}  C(A/J_g)
$
are $G$-invariant subspaces, hence it suffices to prove
the quasiisomorphism
$$
C\Big(\bigoplus_{g \in \mathcal{O}} A^\natural_g \Big)_G \to
\Big(\bigoplus_{g \in \mathcal{O}} C(A/J_g)\Big)_G
$$
for all conjugacy classes $\mathcal{O}$. To that end, choose $g \in \mathcal{O}$
and denote by $C_g = \{h \in G | gh = hg \}$ the centralizer of $g$.
Then $C_g$ acts on $A^\natural_g$; and on the space of
coinvariants $\Big(A/J_g\Big)_{C_g}$ the cyclic operator $t_n$
satisfies $t_n^{n+1} = 1$ therefore we obtain a mixed complex
$C(A^\natural_g)_{C_g}$. Moreover,
we have natural isomorphisms of $G$-modules
$$
 \bigoplus_{g' \in \mathcal{O}} A^\natural_{g'} \simeq Ind_{C_g}^G C(A/J_g);
 \qquad \bigoplus_{g' \in \mathcal{O}} C(A/I_{g'})
 \simeq Ind_{C_g}^G C(A/J_g)
$$
where $Ind_{C_g}^G$ denotes the induction map from $C_g$-modules
to $G$-modules. By Shapiro's Lemma we are reduced to
proving the quasiisomorphism
$$
C(A^\natural_g)_{C_g} \simeq (C(A/J_g))_{C_g};
$$
which would follow once we prove that the natural surjection $A \to A/J_g$
induces a quasiisomorphism
$$
C(A^\natural_g)_{\cyc} \simeq C(A/J_g)
$$
where $\cyc \subset C_g$ is the cyclic subgroup generated by
$g$. Note that $\cyc$ acts trivially on $A/J_g$ (in fact, $J_g$
is generated by elements $a - g(a)$ with $a \in A$), so we don't
have to take coinvariants on the right hand side.

Now we finally use Proposition \ref{mayer} (and the smoothness
assumption which stands behind it). Since the properties (a), (b) and
(c) are preserved by finite intersections of
$G$-invariant affine open subsets, by Mayer-Vietoris argument
we can assume that there exists a point $x \in X$ and a $G_x$-equivariant
map $\varphi: X = Spec\; B \to T_x = Spec \; A$ which is etale
and satisfies (a), (b) and (c) of Proposition \ref{mayer}.
Let $J_B$ and $J_A$ denote the ideal of $g$-fixed points in $B$ and $A$,
respectively. We have two
cases:

\bigskip
\noindent
When $g \notin G_x$, $g$ has no fixed points by (a) of \ref{mayer}, i.e.
$J_B = B$. Then by Theorem 6 of \cite{Lo} $C(B \rtimes \cyc)$ is
quasiisomorphic to $(C(B))_{\cyc}$ (the piece obtained from the
conjugacy class of identity), so the argument of Step 4 applied to
$G = \cyc$ shows that $C(B^\natural_g)$ is quasiisomorphic to zero
(if $g \neq 1$).

\bigskip
\noindent
When $g \in G_x$ we use the etale map $\varphi: X \to T_x$ to reduce to
the case of the flat space $T_x$. Note that until now all maps were defined
in the derived category of $\Lambda$. However,  it
suffices to check that the map $C(B^\natural_g)_{\cyc} \to C(B/J_B)$
defined above is a composition of quasiisomorphisms of complexes
over $k$. Note that the action $b(b_0, \ldots, b_n)
= (b b_0, b_1, \ldots, b_n)$ actually turns both $C(B^\natural_g)$ and $C(B/J_B)$
into complexes of $B$-modules. We will show that they are isomorphic in the
derived category of $B$ (thus, taking coinvariants in $C(B^\natural_g)_{\cyc}$
is only necessary to define an \textit{apriori} mixed complex structure).

\begin{prop} Let $X = Spec \; B$, $x \in X$, $g \in G_x$ and $T_x = Spec\; A$
be as above. If the map $A \to A/J_A$ induces a quasiisomorphism
$C(A^\natural_g) \to C(A/J_A)$ in the derived category of $A$ then the map
$B \to B/J_B$ induces a quasiisomorphism $C(B^\natural_g) \to C(B/J_B)$ in
the derived category of $B$.
\end{prop}
\noindent
\textit{Proof.} This proof is a minor modification of the etale descent
result of \cite{GW}. In fact, if $C(A^\natural_g) \to C(A/J_A)$ is
a quasiisomorphism of complexes of $A$-modules then consider the following
commutative diagram:
$$
\begin{array}{ccc}
 B \otimes_A C(A^\natural_g) & \longrightarrow & B \otimes_A C(A/J_A)
 \simeq B/J_B \otimes_{A/J_A} C(A/J_A) \\
 \downarrow & & \downarrow \\
 C(B^\natural_g) & \longrightarrow & C(B/J_B)
\end{array}
$$
where the two vertical arrows are given by $b \otimes (a_0, \ldots, a_n)
= (ba_0, a_1, \ldots, a_n)$.

The top arrow is a quasiisomorphism since $B$ is flat over $A$. The top left
corner isomorphism holds since $B \simeq B/J_B \otimes_{A/J_A} A$ by part (b)
of Proposition \ref{mayer}. To show that the bottom arrow is a
quasiisomorphism we will show that this property holds for the two vertical arrows.
Moreover, it suffices to prove it for the left arrow only since then one
can set $g =1$, replace the pair $(B, A)$ by $(B/J_B, A/J_A)$, respectively, and
get the proof for the right arrow.

To prove the assertion about $B \otimes_A C(A^\natural_g) \to C(B^\natural_g)$
we borrow some formulas from pp. 513-514 of \cite{BG}. Let $A_{\Delta}$
and $A_g$ denote bimodules over $A \otimes A$ isomorphic to $A$ as abelian
groups, with the module structure given by $(a_0, a_1) \cdot a = a_0 a a_1$
for $a \in A_{\Delta}$ and $(a_0, a_1) \cdot a = a_0 a g^{-1} (a_1)$
for $a \in A_g$. Consider the free resolution $P_\bullet^A$ of $A_{\Delta}$ as
$A \otimes A$-module:
$$
\ldots \to A \otimes A \otimes A \otimes A \stackrel{\partial}\to A \otimes A
\otimes A \stackrel{\partial}\to A \otimes A \stackrel{\Delta}\to  A \to 0
$$
where the $A \otimes A$-module structure on $A^{\otimes (k+2)}$ is given by
$$
(\bar{a}_0, \bar{a}_1) \cdot (a_0, \ldots, a_{k+1}) = (\bar{a}_0 a_0, a_1,
\ldots, a_k, a_{k+1} \bar{a}_1);
$$
the map $\Delta: A \otimes A \to A$ sends $(a_0, a_1)$ to $a_0 a_1$ and
$$
\partial (a_0, \ldots, a_{k+1}) =
\sum_{i = 0}^k (-1)^i (a_0, \ldots, a_i a_{i+1}, \ldots a_{k+1}).
$$
Exactness of $P_\bullet^A$ is proved using the homotopy $s: (a_0, \ldots,
a_{k+1}) \mapsto (1, a_0, \ldots, a_{k+1})$. Let $B_g$ and $P^B_\bullet$ be
the similar objects over $B$.

Now we have a chain of isomorphisms:
$$
B \otimes_A C (A^\natural_g)
\simeq B \otimes_A A_g \otimes_{A \otimes A} P^A_\bullet \simeq
B_g \otimes_{A \otimes A} P^A_\bullet \simeq B_g \otimes_{B \otimes B}
(B \otimes B) \otimes_{A \otimes A} P^A_\bullet.
$$
The first isomorphism follows from the definitions of $P^A_\bullet$, $C(A^\natural_g)$ and
$A_g$. The second uses
$\cyc$-equivariance of $Spec\; B \to Spec\; A$ and the $A \otimes A$-module
structure on $B_g$ which comes from $A \otimes A \to B \otimes B$. Taking
into account that $C(B^\natural_g) \simeq B_g \otimes_{B \otimes B} P^B_\bullet$,
we need to prove that the natural injective map of complexes
$\rho: (B \otimes B) \otimes_{A \otimes A}
P^A_\bullet \to P^B_\bullet$ becomes a quasiisomorphism after applying
$B_g \otimes_{B \otimes B}(\ldots)$.  Explicitly, we have
$$
\begin{array}{cccccccc}
\ldots \to & B \otimes A \otimes B & \to & B \otimes B & \to & B \otimes_A B & \to
& 0 \\
&  \downarrow \rho & & \downarrow \rho & & \downarrow a & \\
\ldots \to & B \otimes B \otimes B & \to & B \otimes B & \to & B & \to & 0
\end{array}
$$
where the upper row is exact since $B \otimes B$ is flat over $A \otimes A$ and the
right vertical arrow $a$ is the natural surjective map $B \otimes_A B \to B$.

Since  $B$ is etale over $A$, the kernel $C$ of $a$ is a $B \otimes B$-module
supported away from the diagonal $X_\Delta \subset X \times X = Spec(B \otimes B)$.
Note that the support of $B_g$ coincides with the graph of the map $g^{-1}: X \to X$.
We now claim that the supports of $B_g$ and $C$ are disjoint. In fact, if $(x_1, x_2)
\in Supp (C) \subset X \times X$ then $x_1 \neq x_2$ but $\varphi(x_1) =
\varphi(x_2) \in T_x$. If $(x_0, g^{-1}(x_0)) \in  Supp(C)$ then $g^{-1}$ stabilizes
$\varphi(x_0) \in T_x$  but does not stabilize $x_0 \in X$, which contradicts
property (c) in Proposition \ref{mayer}. Hence the supports of $B_g$ and $C$ are
disjoint and after tensoring $B_g \otimes_{B \otimes B} $ the two rows
in the above diagram become quasiisomorphic since $B_g \otimes_{B \otimes B} Coker
\rho$ computes $Tor^\bullet_{B \otimes B} (B_g, C) = 0$
(a bit more rigorously, one could first show that the
localization at each maximal prime vanishes - see the last lines on p. 374 of
\cite{GW}).  This finishes the proof of the proposition.
 $\square$

\bigskip
\noindent
Thus, it suffices to prove the quasiisomorphism $C(A^\natural_g) \to C(A/J_A)$ when
$A = k[T_x]$ is a polynomial ring with an action of the cyclic group $\cyc$
induced from its linear action on $T_x$.

\vspace{1cm}

\noindent
\label{step6}
\textit{Step 6: the Linear Case}

\vspace{1cm}
\noindent
As a last step we consider the case of a flat space $V = T_x= Spec\; A$ with a
linear action of the cyclic group $H=\cyc$. Let $V = V_0 \oplus V_1$ where
$V_0 = V^H$ and $V_1$ is the $H$-invariant complement. Note that
$\overline{A}=A/J_A$ is the algebra of regular functions on $V_0$.

Recall that for any variety $Y$, a vector bundle $E$ and its section $s$
one has a Koszul complex
$$
\ldots \to \Lambda^3 E^* \stackrel{\partial}\to \Lambda^2 E^* \stackrel{\partial}
\to E^* \stackrel{\partial} \to \mathcal{O}_Y \to 0
$$
where the differential is given by contraction with $s$. We denote this
Koszul complex by $Kos(Y, E, s)$.  It is well-known that for a regular
section $s$ and affine $Y$, $Kos(Y, E, s)$ is a projective resolution of
$\mathcal{O}_Z$ where $Z$ is the zero scheme of $s$.

Recall from the previous step that $C(A^\natural_g)$ is obtained by
taking a particular projective resolution $P^A_\bullet$ of the
diagonal copy $V_\Delta \subset V \times V$ and tensoring it with
the $A \otimes A$-module corresponding to the graph of $g^{-1}: V \to V$.
Similarly, $C(\overline{A})$ is obtained by taking a particular resolution
$P^{\overline{A}}_\bullet$ of $(V_0)_\Delta \subset V_0 \times V_0$ and
tensoring it with the module of functions on $(V_0)_\Delta$.  As in \cite{BG}
we prove the quasiisomorphism $C(A^\natural_g) \to C(\overline{A})$ by
looking at the Koszul resolutions instead of $P_\bullet$.

Indeed, $V_\Delta \subset V \times V$ is a zero scheme
of a section $s$ of the trivial vector bundle with fiber $V$, given by
$s(v_1, v_2) = v_1 - v_2$; and similarly $(V_0)_\Delta$ is a
zero scheme of the section $s_0 = s|_{V_0 \times V_0}$ which takes
values in the trivial vector bundle with the fiber $V_0$. Then
  we have a commutative diagram
$$
\begin{array}{ccc}
P^A_\bullet \quad & \longrightarrow & P^{\overline{A}}_\bullet \quad \\
\downarrow \alpha_A & & \downarrow\alpha_{\overline{A}} \\
Kos(V \times V, V, s) & \longrightarrow &Kos (V_0 \times V_0, V_0, s_0)
\end{array}
$$
where $\alpha_A$ is an extension of the identity map $A_\Delta \to A_\Delta$
to the projective resolutions, and $\alpha_{\overline{A}}$ is its reduction
modulo $J_A$ (for example, we define $\alpha$ by the
formula on p. 515 of \cite{BG}). Since $\alpha_A$ (resp. $\alpha_{\overline{A}}$)
is a quasiisomorphism by definition of a projective resolution, and remains
one after applying $A_g \otimes_{A \otimes A}$ (resp. $\overline{A}_\Delta
\otimes_{\overline{A} \otimes \overline{A}}$) we only have to check that the
induced map
$$
A_g \otimes_{A \otimes A} Kos(V \times V, V, s) \to \overline{A}_\Delta
\otimes_{\overline{A} \otimes \overline{A}} Kos (V_0 \times V_0, V_0, s_0)
$$
is a quasiisomorphism. But by direct computation
(see \cite{SGA}, Expose VII, Proposition 2.5, for example) one can see that
the left hand side is isomorphic to the left hand side tensored by
$Kos(V_1, V_1, s')$
where $s'$ is the section of a trivial vector bundle
over $V_1$ with fiber $V_1$, given by $s'(v) = v - g^{-1}(v)$. Since $V_1^H = 0$,
$Kos(V_1, v_1, s')$ is quasiisomorphic to $k$. This finishes the proof of Step 6,
and the proof of Proposition \ref{actual}. $\square$

\section{Equivalences of Derived Categories and Cohomology}

\noindent
First we state a result which says that the cohomology of a complex
quasiprojective variety can be recovered from an enhanced version
of its derived category of vector bundles. The parts (a) and (b)
of the next theorem are not stated
explicitly in the papers \cite{K1}, \cite{K2} but easily follow
from their results.

\begin{theorem}
\label{cohom-der}(a) Let $X$ be a quasiprojective variety over the field
of complex numbers. In notation of Section 1, the cyclic homology group
$HC_i(Vect(X), \mathbb{C}[u, u^{-1}])$
is isomorphic to $H^{even}(X, \mathbb{C})$ for $i = 2k$ and
$H^{odd}(X, \mathbb{C})$ for $i = 2k+1$.

(b) If $F: D^b(X) \to D^b(Y)$ is an equivalence of bounded derived
categories of sheaves on two smooth projective varieties $X, Y$ over
a filed $k$ then
$F$ induces an isomorphism $HC_\bullet(Vect(X), W) \to
HC_\bullet(Vect(Y), W)$ for any graded $k[u]$-module $W$
of finite projective dimension. In particular, if $k =\mathbb{C}$ then
$F$ induces an isomorphism of complex cohomology groups.

(c) Let $X, Y$ be as in (b) and assume that a finite group $G$
acts on $X$. If $char\; k$ does not divide $|G|$ and there exists
an equivalence of derived categories $F: D^b(Y) \to D^b_G (X)$
 then $F$ induces an isomorphism
$$
HC_\bullet(Vect(Y), W) \simeq \Big( \bigoplus_{g \in G} HC_\bullet
(Vect(X^g), W) \Big)_G
$$

(d) If $char\; k = 0$ (resp. $k = \mathbb{C}$) and $W = k [u, u^{-1}]$
then in the situation of (b) or (c) one has an isomorphism
of Hodge filtrations (resp. pieces of the Hodge decomposition).
\end{theorem}
\noindent
\textit{Proof.} By Corollary 5.2 of \cite{K2} the mixed complex
$C(X)$ is quasiisomorphic to the mixed complex obtained by
sheafifying the standard mixed complex of an algebra
(see Section 9 of \cite{W1}). This result holds for
any field $k$. If $char\; k = 0$ then by \cite{FT} (in the affine
case) and \cite{W3} (general quasiprojective schemes)
the homology groups $HC_i(Vect(X), k[u, u^{-1}])$
are given by the crystalline cohomology of $X$ (see Theorem 3.4 in \cite{W3}).
If particular, when $k = \mathbb{C}$ we get the Betti cohomology of
$X(\mathbb{C})$. This proves (a).

To prove (b) note that by smoothness the triangulated category obtained from
$(\mathcal{A}c^b(X), \mathcal{C}^b (X))$ is equivalent
to the bounded derived category $D^b(X)$. By
a fundamental theorem of Orlov (see Theorem 2.2,
\cite{Or}) any equivalence $F$ as above is induced by a Fourier-Mukai
transform with respect to some sheaf on $X \times Y$. Thus, any such
equivalence $F$ automatically comes from a morphism of localization pairs
and we can conclude by using the invariance property
(Theorem \ref{invariance}). Note that part (a) gives an isomorphism
of complex cohomology groups only as $\mathbb{Z}/2 \mathbb{Z}$-graded
vector spaces but by the main result of \cite{W3} this can be
refined to an isomorphism of $\mathbb{Z}$-graded vector spaces as
well.

Part (c) follows from Theorem \ref{invariance}; Theorem \ref{main}
and the fact that any equivalence $F$ still comes from a functor
between localization pairs due to 8.1 and 8.2 in \cite{K3}.
Part (d) follows from  \cite{W3}.
$\square$

\bigskip
\noindent
\textbf{Example.}  Let $X -  \to Y$ be an elementary flop of
Bondal-Orlov, see \cite{BO}, Theorem 3.6. By \emph{loc. cit.} any
such flop induces an equivalence
of derived categories $F: D^b(X) \to D^b(Y)$ hence an isomorphism
of cohomology groups by part (b) of the above theorem. Note that the
 motivic integration approach only gives an equality of
Betti numbers or classes in the K-theory of Hodge structures.  This
equivalence
was extended to more general flops in dimensions 3 and 4
by Bridgeland and Namikawa.

\bigskip
\noindent
Next we describe some situations when a derived category of sheaves
on one variety is equivalent to the equivariant derived category of sheaves
(with respect to a finite group action) on another variety. From now
on we assume that $k = \mathbb{C}$. Let $G$ be a finite group with a
unimodular action on a smooth irreducible quasiprojective variety $X$
(i.e. we require that for each $x \in X$ the image of the stabilizer
$G_x$ in $GL(T_x)$ actually belongs to the subgroup $SL(T_x)$). Denote
by $GHilb(X) \subset Hilb^{|G|}(X)$ the scheme of all $G$-invariant
0-dimensional subschemes $Z \subset X$ on multiplicity $|G|$ such that
$H^0(\mathcal{O}_Z)$ is isomorphic to a regular representation of $|G|$.
Assume further that the generic point of $x \in X$ has trivial stabilizer
and let $Y(G, X) \subset GHilb(X)$ be the irreducible component
containing all free $G$-orbits.

\begin{theorem} \label{exmpl-der}  \ \\
\noindent
(1) Assume that $X, G$ and $Y= Y(G, X)$ are as above and one of the following
conditions is satisfied

(a) $\dim X \leq 3$;

(b) $X$ is a complex symplectic variety, $G$ preserves the symplectic
structure and $Y$ is a crepant resolution of $X$;

(c) $X$ is the $n$-th cartesian power of a smooth quasiprojective surface
$S$; $G = \Xi_n$ (the symmetric group) with the natural permutation
action on $X$.

\medskip
\noindent
Then there exists an equivalence $F: D^b(Y) \to D^b_G(X)$ coming from
a morphism of localization pairs.

\medskip
\noindent
(2) Assume that $\dim X \leq 3$; $X$ is projective and $X'$ is another
smooth projective variety with an action of a finite group $G'$. Suppose
that $X/G$ and $X'/G'$ are Gorenstein and there exists a common resolution
of singularities $\pi: Z \to X/G$,
$\pi': X'/G'$ with $\pi^* K_X \simeq (\pi')^* K_{X'}$. Then there
exists an equivalence of categories $D^b_G(X) \simeq D^b_{G'}(X')$.
In particular, if $Z = X'$ and $G'= \{1\}$ (i.e. $Z$ is a crepant
resolution) there exists an equivalence $D_G^b(X) \simeq D^b(Z)$.
\end{theorem}
\noindent
\textit{Proof.}  Parts 1a, 1b correspond to  Theorem 1.2 and Corollary 1.3
in \cite{BKR}, respectively. To prove 1c, note that for $S = \mathbb{C}^2$
by a result of Haiman, cf. \cite{Ha} Theorem 5.1, the variety $Y$ is
isomorphic to the Hilbert scheme $Hilb^n (S)$. Then the morphism $Y \to
X/\Xi_n$ is semismall, see \cite{GS}, and the assertion follows from a
general criterion of Theorem 1.1 in \cite{BKR}. In general the equality
$Hilb^n(S) = Y$ is derived from the $\mathbb{C}^2$ case by considering
completions of local rings at points $s \in S$ which are the same as for
$\mathbb{C}^2$ since $S$ is smooth.

Part 2 is a simplified version of Theorem 1.7 due to Kawamata, cf. \cite{Ka}.
$\square$.

\begin{corr} \label{last}In the first (resp. the second) case of the above theorem one
has an isomorphism
of the cyclic homology groups
$$
 HC_\bullet(Vect_G(X), W) \simeq HC_\bullet(Vect(Y), W), \quad (resp. \quad
 HC_\bullet(Vect_G(X), W) \simeq HC_\bullet(Vect_{G'}(X'), W))
$$
for all graded $k[u]$-modules $W$ of finite projective dimension.
For $W = \mathbb{C}[u, u^{-1}]$ this reduces to isomorphisms
$$
H^*_{orb} (X, \mathbb{C}) \simeq H^*(Y, \mathbb{C}); \qquad
H^*_{orb} (X, \mathbb{C}) \simeq H^*_{orb}(X', \mathbb{C})
$$
of (orbifold) cohomology groups with their Hodge filtratons.
\end{corr}
\noindent
\textit{Proof.} By inspecting the proofs in \cite{BKR} and and \cite{Ka}
one can see that all the derived equivalences above are given
by functors between localization pairs. Hence the first two assertions
follow from the invariance property (Theorem \ref{invariance}).
To obtain the last pair of equalities on applies Theorem \ref{main} and
Theorem \ref{cohom-der} (a). $\square$

\medskip
\noindent
Finally we note that by a recent conjecture of Kawamata, the derived
equivalences  \ref{exmpl-der} should be a part of more general statement.
Here we formulate a part of  Conjecture 1.2 in \cite{Ka} in a slightly
generalized form:

\begin{conj} (Categorical Resolution Conjecture)
Let $\mathcal{X}$, $\mathcal{Y}$ be  smooth Deligne-Mumford
stacks with Gorenstein moduli spaces $X, Y$. Suppose that
there exist
birational maps $f:Z \to X$ and $g: Z \to Y$ such that
$f^*(K_X) = g^*(K_Y)$. Then the derived categories $D^b(\mathcal{X})$ and
$D^b(\mathcal{Y})$ of sheaves in the etale topology, are equivalent.
\end{conj}

\section{Concluding Remarks}

The isomorphisms of Theorem \ref{main} are additive counterparts
of a $K$-theoretic statement in \cite{V}. However, in \textit{loc. cit.}
one also has a statement concerning $K'$-theory of singular varieties.
This motivates the following conjecture:

\begin{conj} Let $X$ be a quasiprojective scheme over a field
$k$, $G$ a finite group acting on $X$ and assume that $char\; k$ does not
divide $|G|$. Let $Coh(X)$ be the exact category of coherent sheaves on $X$
and $Coh_G(X)$ the exact category  of $G$-equivariant sheaves.
Then for any $W$ of finite projective dimension over $k[u]$
there exists an isomorphism
$$
\phi_X:   \Big( \bigoplus_{g \in G}
HC_\bullet(Coh(X^g), W) \Big)_G   \to HC_\bullet(Coh_G(X), W)
$$
which is functorial with respect to the (derived) pushforwards under
$G$-equivariant proper maps.
\end{conj}

We expect that for smooth $X$ and $W = k[u, u^{-1}$ (periodic cyclic homology
case) this isomorphism satisfies  the multiplicative part of the Cohomological
Crepant Resolution Conjecture, see
\cite{Ru}.

We have seen in the proof of Theorem \ref{main} that the isomorphism
$\psi_X$ is defined using pullbacks to the fixed point sets. We expect
that, as in \cite{V}, the isomorphism $\phi_X$ can be defined using
direct images with respect to the closed embeddings $X^g \to X$. Note
that $\phi_X$ and $\psi_X$ are not expected to be mutually inverse, but
their composition should be invertible similarly to
Lemma 4.2 and Lemma 4.3 in \cite{V}. For singular $X$ the
cyclic homology of $Coh(X)$ may be different from the sheafified
cyclic homology of the algebras $\mathcal{O}(U)$, $U \subset X$; therefore
the general strategy of proof has to be completely different.

\bigskip
\noindent
Finally we state another conjecture aimed at a better understanding of the
orbifold product on orbifold cohomology, see Section 3 \cite{Ru}.

Let $\mathcal{X}$ be a smooth Deligne-Mumford stack and $\pi: U \to
\mathcal{X}$ be an etale cover by  a scheme. Consider $Z = U \times_{\pi} U
\subset U \times U$. If $\pi_i: U \times U \times U \to U \times U$
is the projection omitting the $i$-th factor, $i = 1, 2, 3$, then
$\pi_2(\pi_3^{-1} (Z) \cap \pi_1^{-1}(Z)) = Z$, i.e. $Z$ is an
idempotent correspondence from $U$ to itself.

In particular, we have a (finite) morphism $m: Z \times_U Z \to Z$
which allows to define for any two objects $\mathcal{F}$,
$\mathcal{G} \in D^b(Z)$ a third object
$$
\mathcal{F} * \mathcal{G} := m_* \big(p_1^*\mathcal{F} \otimes p_2^*\mathcal{G})
$$
where $p_1, p_2$ are the projections of $Z \times_U Z$ to $Z$, and all operations
are understood in the derived sense.

The $*$-product in general is not commutative, but it is
associative (more precisely, there is a canonical isomorphism between
$\mathcal{F}*(\mathcal{G} * \mathcal{H})$ and $(\mathcal{F} * \mathcal{G})
* \mathcal{H}$). When $\mathcal{X}$ is the quotient stack $[X/G]$ we
can take $U = X$ and then $Z = \coprod_{g \in G} X^g$ while
$m: Z \times_U Z \to Z$ is given by pointwise group product within
stabilizers.

\begin{conj} When $\mathcal{X} = [X/G]$ and $U = X$ the above
product $(\mathcal{F}, \mathcal{G}) \mapsto \mathcal{F} * \mathcal{G}$
induces the product of the periodic cyclic homology of $Z = \coprod_{g \in G} X^g$
which coincides with the product on $A(X, G) = H^*(Z, \mathbb{C})$ described
in Section 3 of \cite{Ru}. Thus, taking the (co)invariants with respect to
$G$-action gives the orbifold product structure on $H^*_{orb}(X/G, \mathbb{C})$.
\end{conj}

Note that the only place where smoothness of  $\mathcal{X}$ is important
is the relation between cyclic homology of coherent sheaves and the
orbifold cohomology groups. However, we could \textit{define} orbifold
cohomology for singular Deligne-Mumford stacks via cyclic homology of
coherent sheaves; such a definition, perhaps, would give invariants which
behave better than those defined via the usual cohomology.

\bigskip

Department of Mathematics

253-37 Caltech

Pasadena, CA 91125, USA

\smallskip

email: baranovs@caltech.edu

\end{document}